\documentclass[12pt]{article}
\def\q{\hfill\rule{1ex}{1ex}}
\def\0{\emptyset}

\def\p{{\bf Proof.} \quad}
\def\q{\hfill\rule{1ex}{1ex}}

\def\n{\noindent}

\usepackage{graphicx}

\begin{document}
\title{\bf Edge-bipancyclicity of bubble-sort star graphs}
\author{{\small\bf
Jia Guo$^1$}\thanks{email: guojia199011@163.com}\quad\quad{\small\bf
Mei Lu$^2$}\thanks{email: lumei@mail.tsinghua.edu.cn}\\
{\small $^1$College of Science, Northwest A$\&$F University, Yangling, Shaanxi 712100, PR China}\\
{\small $^2$Department of Mathematical Sciences, Tsinghua University, Beijing
100084, PR China}
}

\date{}
\maketitle\baselineskip 15.5pt

\begin{abstract}
\baselineskip=0.5cm The interconnection network considered in this
paper is the bubble-sort star graph. The $n$-dimensional bubble-sort star graph $BS_n$ is a bipartite and $(2n-
3)$-regular graph of order $n!$. A bipartite graph $G$ is
edge-bipancyclic if each edge of $G$ lies on a cycle of all even length $l$ with $4\leq l\leq |V(G)|$.
In this paper, we show that the $n$-dimensional bubble-sort star graph $BS_n$ is edge-bipancyclic for $n\ge 3$ and for each even length $l$ with $4\leq l\leq n!$,
every edge of $BS_n$ lies on at least four different cycles of length $l$. 

\vskip 0.3cm

{\bf AMS} classification: 05C50  \vskip 0.3cm

{\bf Keywords:} bubble-sort star graph, edge-bipancyclicity

\end{abstract}

\vskip.3cm

\n{\large\bf 1.\quad Introduction}

Network topology is a crucial factor for the interconnection
networks since it determines the performance of the networks. Many
interconnection network topologies have been proposed in the
literature for the purpose of connecting hundreds or thousands of
processing elements. In these interconnection network topologies,
linear arrays and rings are two of the most fundamental structures.
One of the important issues in parallel processing is to embed
linear arrays and rings into interconnection networks (see
\cite{Latifi} and \cite{Sengupta}). Paths and cycles are popular
interconnection networks owing to their simple structures and low
degree. Moreover, many parallel algorithms have been devised on them
\cite{Lin, Hallaron}. Many literatures have addressed how to embed
cycles and paths into various interconnection networks
\cite{Auletta, Chen, Fan, Fang, Liu}.

An interconnection network is usually represented by an undirected
simple graph where vertices represent processors and edges represent
links between processors. In the rest of the paper, we will use
Bondy and Murty \cite{Bondy1} for terminology and notation not
defined here and only consider simple undirected graphs.

Let $G=(V,E)$ be a graph. $G$ is called {\em $k$-regular} if every vertex in $G$ has exactly $k$ neighbors. $P_k=v_1v_2\cdots v_k$
on $k$ distinct vertices $v_1,v_2,\cdots,v_k$ in  $G$
is a {\em $k$-path} if $(v_i,v_{i+1})$ is an edge in $G$ for every $i=1,2,\cdots,k-1$.
For a $k$-path $P_k=v_1v_2\cdots v_k$, if $(v_1,v_k)$ is an edge, then $C_k=v_1v_2\cdots v_k v_1$ is a {\em $k$-cycle} in $G$ and the {\em length} of  $C_k$ is $k$.
$G$ is {\em Hamiltonian} if it contains a cycle passing through
every vertex exactly once and the cycle is called a {\em Hamiltonian cycle}.
We say
that $G$ of order $n$ is {\em pancyclic} if it contains a cycle of every possible length
between 3 and $n$.
$G$ is {\em
vertex-pancyclic} (resp.
{\em
edge-pancyclic}) if every vertex (resp. edge) of $G$ lies on a cycle of all length
$l$ with $3\le l\le |V(G)|$.
$G$ is bipartite if the vertex set of $G$ can be partitioned
into two vertex subsets $V_1$ and $V_2$ with $V(G)=V_1\cup V_2$, $V_1\cap V_2=\emptyset$ such that each edge of $G$
joins one vertex in $V_1$ and the other in $V_2$. It is well known
that any bipartite graph contains no odd cycles.
A bipartite graph $G$ of order $n$ is {\em bipancyclic} if it contains a cycle of every possible even length
between 4 and $n$.
A bipartite graph $G$ is {\em
vertex-bipancyclic} (resp.
{\em edge-bipancyclic}) if every vertex (resp. edge) of $G$ lies on a cycle of all even length
$l$ with $4\le l\le |V(G)|$. It is obvious that every edge-pancyclic (resp. edge-bipancyclic) graph is vertex-pancyclic (resp. vertex-bipancyclic) graph and every vertex-pancyclic (resp. vertex-bipancyclic) graph is pancyclic (resp. bipancyclic) graph.
The edge-bipancyclicity and the bipancyclicity of different interconnection networks are widely studied. For
example, see \cite{Kikuchi, Kuo, Li1, Li2, Li, Stewart}.

The bubble-sort star graphs \cite{Chou}, which belong to
the class of Cayley graphs, have been attractive alternative to the
hypercubes. It gains many nice properties, such as high degree of
regularity and symmetry. In particular, the $n$-dimensional
bubble-sort star graph $BS_n$ has $n!$ vertices, and is $(2n-3)$-regular
and vertex transitive. But it is not edge transitive. $BS_n$ has received considerable attention in recent years. For example,
Cai et al. \cite{Cai} investigated the fault-tolerant maximal local-connectivity of $BS_n$. Guo et al. \cite{Guo} gave the conditional diagnosability of $BS_n$.
Gu et al. \cite{Gu} determined the pessimistic diagnosability of $BS_n$.
Wang et al. \cite{Wang1, Wang2} studied the 2-extra connectivity (resp. 2-good-neighbor connectivity) and the 2-extra diagnosability (resp. 2-good-neighbor diagnosability) of $BS_n$.

This paper deals with the edge-bipancyclicity of bubble-sort star graphs. In
\cite{Chou}, Chou et al. showed that the bubble-sort star graph
is Hamiltonian.
In this paper, we will show that  $n$-dimensional
bubble-sort star graph $BS_n$ is
edge-bipancyclic, vertex-bipancyclic and bipancyclic for $n\ge 3$ and for each even length $l$ with $4\leq l\leq n!$,
every edge of $BS_n$ lies on at least four different cycles of length $l$. \vskip.3cm

\vskip.3cm

\n{\large\bf 2.\quad Preliminaries}

\vskip.2cm

In this section, we first review bubble-sort star graphs and give some
lemmas which will be used in the following proof.

\vskip.2cm
{\bf Definition 2.1 \cite{Chou}} The $n$-dimensional bubble-sort star graph $BS_n$ has
vertex set that consists of all $n!$ permutations on $\{1,2,\cdots,
n\}$. A permutation $x$ on $\{1,2,\cdots, n\}$ is denoted as $x=
x_1x_2\cdots x_n$. A vertex $x=x_1x_2\cdots x_n\in V(BS_n)$ is adjacent to a vertex $y=y_1
y_2\cdots y_n\in V(BS_n)$ if and only if there exists an integer $i$ with $2\leq i\leq n$ such that $y_i=x_{i-1}$, $y_{i-1}=x_i$ and $x_j=y_j$ for every $j\in\{1,2,\cdots,n\}-\{i-1,i\}$ or $y_1=x_i$, $y_i=x_1$ and $x_j=y_j$ for every $j\in\{2,\cdots,n\}-\{i\}$.

\vskip.2cm

By Definition 2.1, $BS_n$ is a bipartite graph that has $n!$
vertices, each of which is a permutation on $\{1,2,\cdots, n\}$ and
each vertex has degree $2n-3$. Fig. 1 shows $BS_2$, $BS_3$, and $BS_4$,
respectively.

\begin{figure}[ht]
   \begin{center}
    \includegraphics[scale=1]{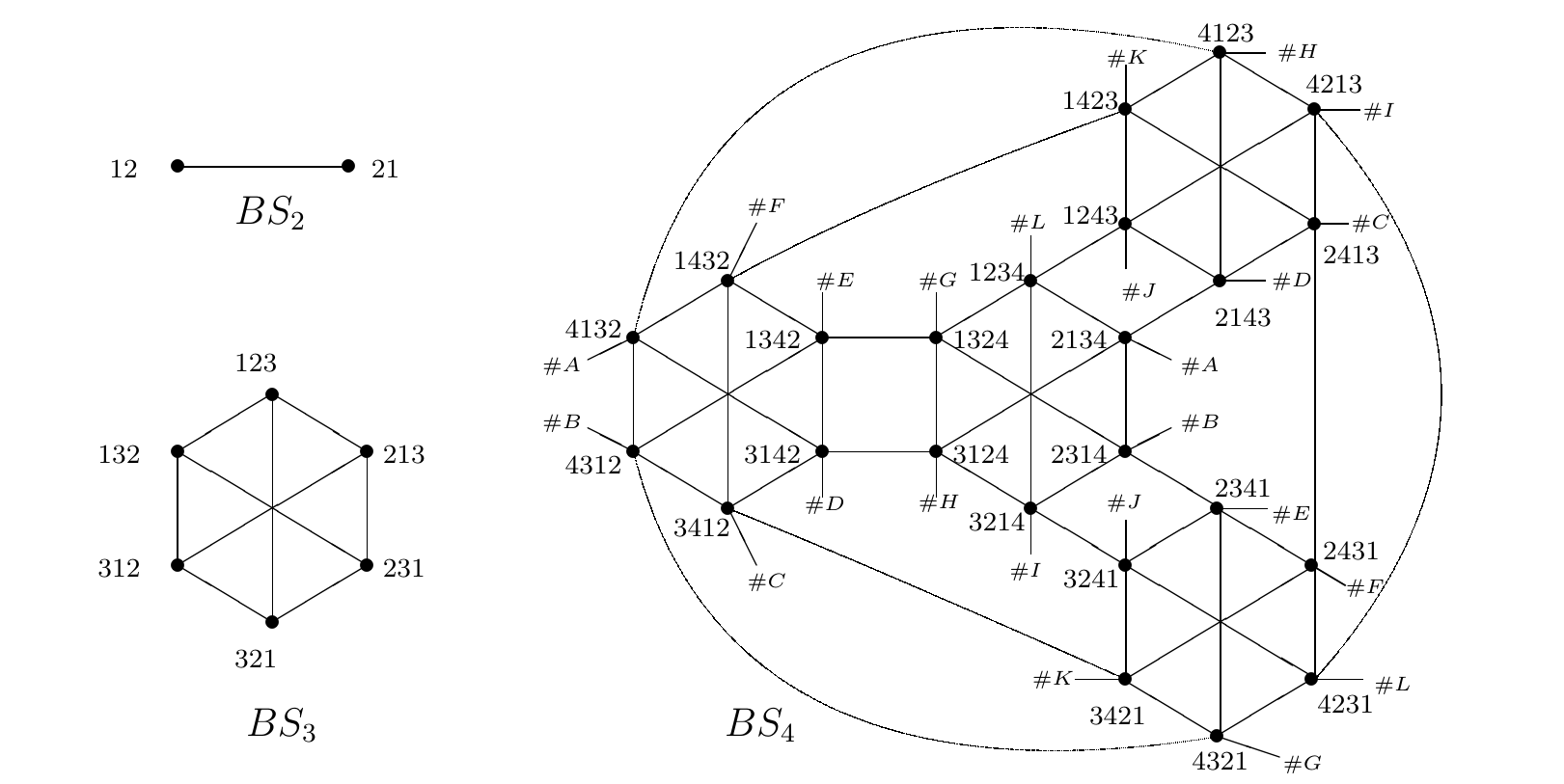}
    \end{center}
    \vspace{-0.5cm}\caption{\label{F2-5} Illustration of $BS_n$ for $n=2,3,4.$}
\end{figure}

\vskip.3cm

Let $x=x_1x_2\cdots x_n,~y=y_1y_2\cdots y_n\in V(BS_n)$,
we use ``$\circ$" to denote an operation such that $x=y\circ(i,j)$ if and only if $x_i=y_j$, $x_j=y_i$ and $x_k=y_k$ for every $k\in\{1,2,\cdots,n\}-\{i,j\}$. Then $(x,y)\in E(BS_n)$ if and only if $y=x\circ(1,i)$ or $y=x\circ(i-1,i)$ for some $i\in\{2,3,\cdots, n\}$. Let $x^+=x\circ(1,n)$ and $x^-=x\circ(n-1,n)$ for simplicity.
For an integer $i$ with $1\le i\le n$, we use $BS_n(i)$ to denote the induced subgraph of $BS_n$ by
the set of vertices $\{x=x_1x_2\cdots x_n~|~x_n=i\}$. By Definition 2.1, $BS_n(i)\cong BS_{n-1}$ for every $i=1,2,\cdots,n$.
Let $e=(x,y)\in E(BS_n(i))$ and $x'\in\{x^+,x^-\}, y'\in\{y^+,y^-\}$. If $e'=(x',y')\in E(BS_n(j))$ for some $j\neq i$, then $e'$ is called a {\em coupled pair-edge} of $e$.

Now we give some properties about $BS_n$.

\vskip.2cm

{\bf Lemma 2.2  \cite{Chou}} {\em For $n\ge 3$, the bubble-sort star
graph $BS_n$ is Hamiltonian.}
\vskip.2cm

{\bf Lemma 2.3 } {\em Let $n,~k,~m$ be three integers with $n\geq4$, $1\leq k,m \leq n$ and $m\neq k$.
Let $H$ be an arbitrary Hamiltonian cycle in $BS_n(k)$. Then for any vertex $u=u_1u_2\cdots u_{n-2}mk$, there is a vertex $v=v_1v_2\cdots v_{n-1}k$ such that
$e=(u,v)$ has a coupled pair-edge $e'\in E(BS_n(m))$, where $(u,v)\in E(H)$.}
\vskip.2cm

\p If $v_{n-1}=m$, then $e'=(u^-,v^-)$ is a coupled pair-edge of $e$ with $e'\in E(BS_n(m))$. Now we suppose that $v_{n-1}\neq m$. Thus the two neighbors of $u$ in $H$ is $u\circ(1,n-1)$ and $u\circ(n-2,n-1)$. In this case, set $v=u\circ(1,n-1)$. Then $e'=(u^-,v^+)$ is a coupled pair-edge of $e$ with $e'\in E(BS_n(m))$.\q

\vskip.3cm

\n{\large\bf 3.\quad Edge-bipancyclicity of $BS_n$}

\vskip.2cm

We first consider the edge-bipancyclicity of $BS_3$.

\vskip.2cm

{\bf Lemma 3.1 } {\em $BS_3$ is edge-bipancyclic.
Furthermore,
every edge of $BS_3$ lies on at least four different cycles of length 4 and 6, respectively.}

\vskip.2cm

\p Let $e$ be an arbitrary edge in $BS_3$. Since $BS_3$ is vertex transitive, without loss of generality, let $123$ be an end vertex of $e$. Then we will consider the following two cases.

{\bf Case 1.} {\em $e=(123,132)$ or $e=(123,213)$.}

Let $e=(123,132)$. The four different 4-cycles containing $e$ are showed in Fig. 2 and the four different 6-cycles containing $e$ are showed in Fig. 4 (a), (b), (c) and (d).

\begin{figure}[ht]
   \begin{center}
    \includegraphics[scale=0.65]{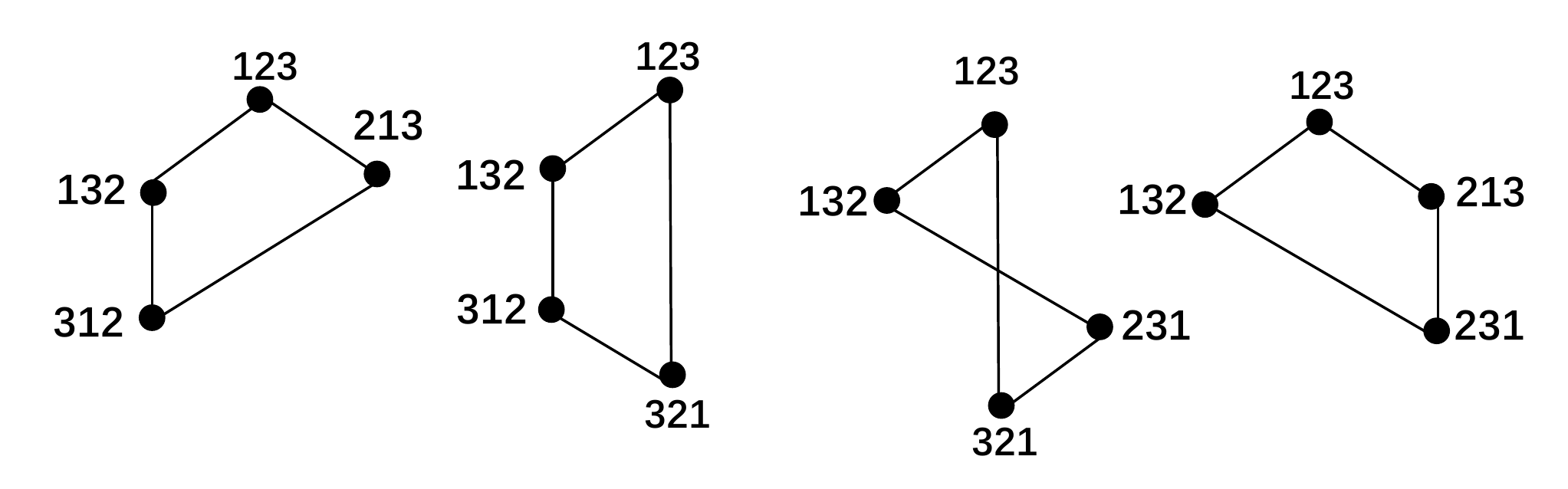}
    \end{center}
    \vspace{-0.5cm}\caption{\label{F2-5} 4-cycles containing $(123,132)$ in $BS_3$.}
\end{figure}

If $e=(123,213)$, the cycles of length 4 and 6 containing $e$ can be constructed similarly to the above case.

\begin{figure}[ht]
   \begin{center}
    \includegraphics[scale=0.65]{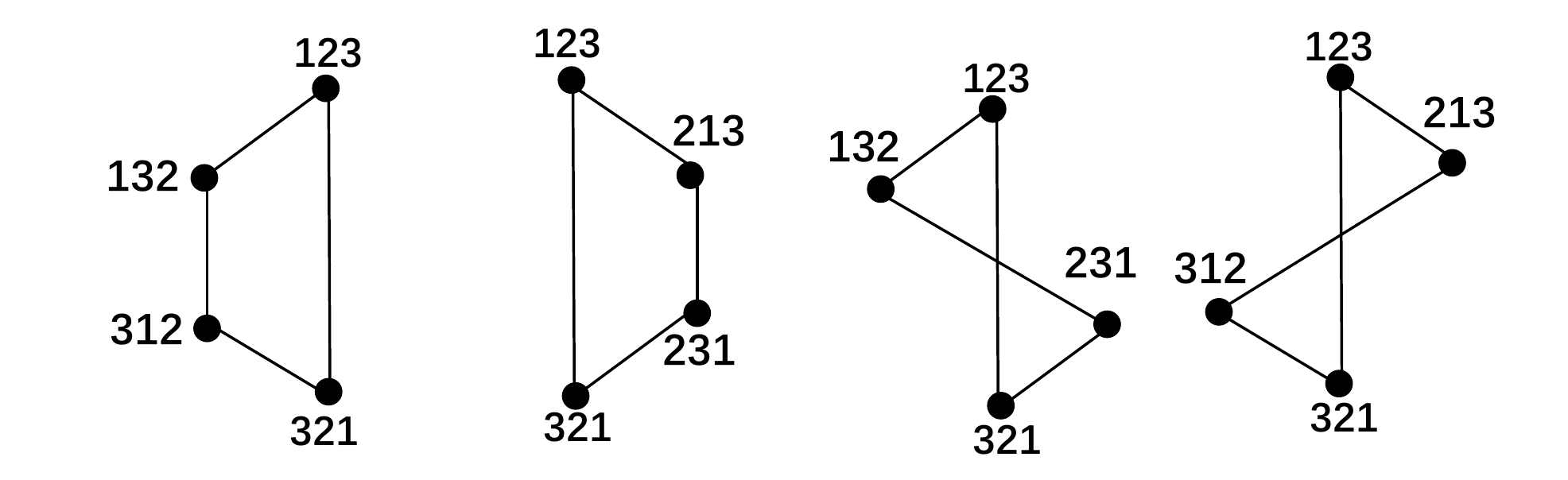}
    \end{center}
    \vspace{-0.5cm}\caption{\label{F2-5} 4-cycles containing $(123,321)$ in $BS_3$.}
\end{figure}

\begin{figure}[ht]
   \begin{center}
    \includegraphics[scale=0.65]{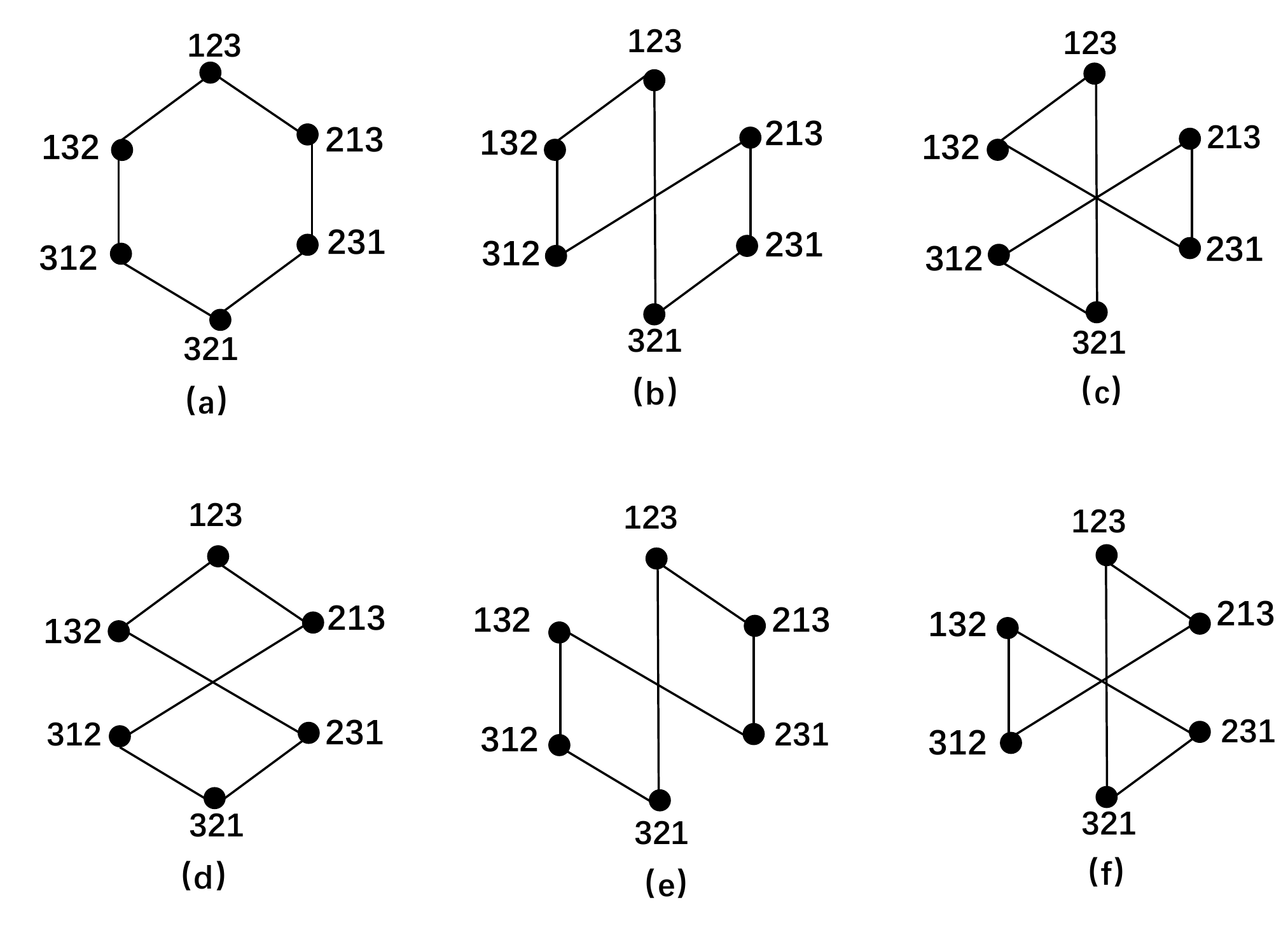}
    \end{center}
    \vspace{-0.5cm}\caption{\label{F2-5} 6-cycles in $BS_3$.}
\end{figure}

{\bf Case 2.} {\em $e=(123,321)$.}

The four different $4$-cycles containing $e$ are showed in Fig. 3 and the four different 6-cycles containing $e$ are showed in Fig.4 (b), (c), (e) and (f).

Thus we get that $BS_3$ is edge-bipancyclic and every edge of $BS_3$ lies on at least four different cycles of length 4 and 6, respectively.\q

\vskip.2cm

Now we show that $BS_4$ is edge-bipancyclic. Lemmas 3.1 and 3.2 are the basis of
induction for proving $BS_n$ to be edge-bipancyclic.
\vskip.2cm

{\bf Lemma 3.2} {\em $BS_4$ is edge-bipancyclic.
Furthermore,
every edge of $BS_4$ lies on at least four different cycles of every even length $l$ with $4\leq l\leq24$.}
\vskip.2cm

\p Let $(u,v)$ be an arbitrary edge of $BS_4$. Since $BS_4$ is
vertex transitive, we can assume that $u=1234$. Then $u\in
V(BS_4(4))$. Now we will consider the following three cases.

{\bf Case 1.} {\em $v\in V(BS_4(4))$.}

Note that $BS_4(4)\cong BS_3$. By Lemma 3.1, $(u,v)$
lies on at least four different cycles of length $4$ and $6$ in $BS_4(4)$, respectively. Now we will construct cycles of every even length between $8$ and $24$ containing  $(u,v)$ by the following two subcases.

{\bf Subcase 1.1.} {\em $v=1324$ or $v=2134$.}

Let $v=1324$.
The four different $8$-cycles containing $(u,v)$
in $BS_4$ are listed in Table 1.

\begin{center}
\setlength{\tabcolsep}{3.5pt}\renewcommand\arraystretch{0.8}
\begin{tabular}{cllllllllllllll}\hline

\multicolumn{1}{c}{$cycle$}&\multicolumn{11}{l}{\scriptsize }
\\\hline

 $C_8^1$&1234&1324&3124&3214&2314&2134&2143&1243&1234&&
\\

$C_8^2$&1234&1324&3124&3214&2314&4312&4132&2134&1234
\\

$C_8^3$&1234&1324&3124&3214&3241&2341&2314&2134&1234
\\

$C_8^4$&1234&1324&3124&4123&4213&3214&2314&2134&1234

\\\hline

\end{tabular}

\vskip 0.2cm

Table 1.  $8$-cycles containing  $(1234,1324)$ in $BS_4$.
\end{center}

Note that $(2143,1243)\in E(BS_4(3))$ is an edge of $C_8^1$ in Table 1. Since $BS_4(3)\cong BS_3$, by Lemma 3.1, (2143,1243) can be contained in four different $t$-cycles in $BS_4(3)$ for $t=4,6$. Let $C^1$ be any $t$-cycle in $BS_4(3)$ containing (2143,1243), then the edge set $E(C_8^1)\cup E(C^1)-\{(2143, 1243)\}$ forms a $(6+t)$-cycle containing (1234,1324). Thus (1234,1324) lies on
at least four different cycles of length 10 and 12, respectively.

Let $C_{12}=(1234)(1324)(3124)(3214)(2314)(2134)(2143)(2413)(4213)(4123)(1423)\\(1243)(1234)$ be a $12$-cycle containing the edge (1234,1324). Let $C_6^1=(1234)(1324)\\(3124)(3214)(2314)(2134)(1234)$, $C_6^2=(2143)(2413)(4213)(4123)(1423)(1243)(2143)$. Thus $C_6^1$ (resp. $C_6^2$) is a Hamiltonian cycle in $BS_4(4)$ (resp. $BS_4(3)$). By Lemma 2.3, there exists an edge $e=(x,y)\not\in\{(1234,2134),(1243,2143)\}$ in $C_6^1$ (resp. $C_6^2$) that has a coupled pair-edge $e'$ such that $e'=(x',y')\in E(BS_n(k))$, where $k=1,2$. Hence the edge set $E(C_{12})\cup \{e',(x,x'),(y,y')\}-\{e\}$ forms a $14$-cycle containing (1234,1324). Thus  (1234,1324) lies on
at least four different 14-cycles.

Let $C_{14}=(1234)(1324)(3124)(3214)(2314)(2134)(2143)(2413)(4213)(4123)(4132)\\(1432)(1423)(1243)(1234)$ be a $14$-cycle containing the edge (1234,1324). Note that $(4132,1432)\in E(BS_4(2))$ is an edge of $C_{14}$. Since $BS_4(2)\cong BS_3$, by Lemma 3.1, $(4132,1432)$ can be contained in four different $t$-cycles in $BS_4(2)$, where $t=4,6$. Let $C^2$ be any $t$-cycle in $BS_4(2)$ containing $(4132,1432)$, then the edge set $E(C_{14})\cup E(C^2)-\{(4132,1432)\}$ forms a $(12+t)$-cycle containing (1234,1324). Thus (1234,1324) lies on
at least four different cycles of length 16 and 18, respectively.

Let $C_{18}=(1234)(1324)(3124)(3214)(2314)(2134)(2143)(2413)(4213)(4123)(4132)\\(4312)(3412)(3142)(1342)
(1432)(1423)(1243)(1234)$ be a $18$-cycle containing the edge (1234,1324). By Lemma 2.3, there exists an edge $e=(x,y)\not\in\{(1234,2134),(1243,\\2143),(1423,4123)\}$ in $C_6^1$ (resp. $C_6^2$) that has a coupled pair-edge $e'$ such that $e'=(x',y')\in E(BS_n(1))$.
Let $C_6^3=(1432)(4132)(4312)(3412)(3142)(1342)(1432)$. Thus $C_6^3$ is a Hamiltonian cycle in $BS_4(2)$.
We have that  $e=(x,y)=(1342,1432)$ (resp. (4312, 3412)) in $C_6^3$ has a coupled pair-edge $e'=(x',y')=(2341, 2431)$ (resp. (4321, 3421)) in $BS_n(1)$. Hence the edge set $E(C_{18})\cup \{e',(x,x'), (y,y')\}-\{e\}$ forms a $20$-cycle containing (1234,1324). Thus  (1234,1324) lies on
at least four different 20-cycles.

Let $C_{20}=(1234)(1324)(3124)(3214)(2314)(2134)(2143)(2413)(4213)(4123)(4132)\\(4312)(4321)(3421)(3412)(3142)(1342)
(1432)(1423)(1243)(1234)$ be a $20$-cycle containing the edge (1234,1324).
Note that $(4321, 3421)\in E(BS_4(1))$ is an edge of $C_{20}$. Since $BS_4(1)\cong BS_3$, by Lemma 3.1, $(4321, 3421)$ can be contained in four different $t$-cycles in $BS_4(1)$ for $t=4,6$. Let $C^3$ be any $t$-cycle in $BS_4(1)$ containing $(4321, 3421)$, then the edge set $E(C_{20})\cup E(C^3)-\{(4321, 3421)\}$ forms a $(18+t)$-cycle containing (1234,1324). Thus (1234,1324) lies on
at least four different cycles of length 22 and 24, respectively.

If $v=2134$, the cycles of every even length between $8$ and $24$ containing (1234,2134) can be constructed similarly to the case $v=1324$.

{\bf Subcase 1.2.}  {\em $v=3214$.}

The four different $8$-cycles containing  $(u,v)$
in $BS_4$ are listed in Table 2.

\begin{center}
\setlength{\tabcolsep}{3.5pt}\renewcommand\arraystretch{0.8}
\begin{tabular}{cllllllllllllll}\hline

\multicolumn{1}{c}{$cycle$}&\multicolumn{11}{l}{\scriptsize }
\\\hline
$C_8^{1*}$&1234&3214&3124&1324&2314&2134&2143&1243&1234
\\

 $C_8^{2*}$&1234&3214&3124&3142&1342&1324&2314&2134&1234&&
\\

$C_8^{3*}$&1234&3214&3124&1324&2314&4312&4132&2134&1234
\\

$C_8^{4*}$&1234&3214&4213&4123&3124&1324&2314&2134&1234

\\\hline

\end{tabular}

\vskip 0.2cm

Table 2.  $8$-cycles containing $(1234,3214)$ in $BS_4$.
\end{center}

Note that $(2143, 1243)\in E(BS_4(3))$ is an edge of $C_8^{1*}$ in Table 2. Since $BS_4(3)\cong BS_3$, by Lemma 3.1, $(2143, 1243)$ can be contained in four different $t$-cycles in $BS_4(3)$ for $t=4,6$. Let $C^{1*}$ be any $t$-cycle in $BS_4(3)$ containing $(2143, 1243)$, then the edge set $E(C_8^{1*})\cup E(C^{1*})-\{(2143, 1243)\}$ forms a $(6+t)$-cycle containing (1234,3214). Thus  (1234,3214) lies on
at least four different cycles of length 10 and 12, respectively.

The cycles of every even length between $14$ and $24$ containing (1234,3214) can be constructed  by the same argument as that of Subcase 1.1.

{\bf Case 2.} {\em $v=u^-=1243$.}

The four different $4$-cycles containing  $(u,v)$
in $BS_4$ are listed in Table 3.

\begin{center}
\setlength{\tabcolsep}{3.5pt}\renewcommand\arraystretch{0.8}
\begin{tabular}{cllllllllllll}\hline

\multicolumn{1}{c}{$cycle$}&\multicolumn{11}{l}{\scriptsize }
\\\hline
$C_4^{1}$&1234&1243&2143&2134&1234
\\

$C_4^{2}$&1234&1243&3241&3214&1234&&&&&&

\\

$C_4^{3}$&1234&1243&4213&3214&1234
\\

$C_4^{4}$&1234&1243&3241&4231&1234

\\ \hline

\end{tabular}

\vskip 0.2cm

Table 3.  $4$-cycles containing  $(1234,1243)$ in $BS_4$.
\end{center}

Note that $(2134, 1234)\in E(BS_4(4))$ is an edge of $C_4^{1}$ in Table 3. Since $BS_4(4)\cong BS_3$, by Lemma 3.1, $(2134, 1234)$ can be contained in four different $t$-cycles in $BS_4(4)$ for $t=4,6$. Let $C^{1'}$ be any $t$-cycle in $BS_4(4)$ containing $(2134, 1234)$, then the edge set $E(C_4^{1})\cup E(C^{1'})-\{(2134, 1234)\}$ forms a $(2+t)$-cycle containing (1234,1243). Thus (1234,1243) lies on
at least four different cycles of length 6 and 8, respectively.

Let $C_8^1$ be the 8-cycle containing (1234,1243) in Table 1. Then the cycles of every even length between $10$ and $24$ containing (1234,1243) can be constructed by the same argument as that of Subcase 1.1.

{\bf Case 3.} {\em $v=u^+=4231$.}

The four different $4$-cycles containing $(u,v)$
in $BS_4$ are listed in Table 4.

Note that $(3214, 1234)\in E(BS_4(4))$ is an edge of $C_4^{1*}$ in Table 4. Since $BS_4(4)\cong BS_3$, by Lemma 3.1, $(3214, 1234)$ can be contained in four different $t$-cycles in $BS_4(4)$, where $t=4,6$. Let $C^{1''}$ be any $t$-cycle in $BS_4(4)$ containing $(3214, 1234)$, then the edge set $E(C_4^{1*})\cup E(C^{1''})-\{(3214, 1234)\}$ forms a $(2+t)$-cycle containing (1234,4231). Thus (1234,4231) lies on
at least four different cycles of length 6 and 8, respectively.

\begin{center}
\setlength{\tabcolsep}{3.5pt}\renewcommand\arraystretch{0.8}
\begin{tabular}{cllllllllllll}\hline

\multicolumn{1}{c}{$cycle$}&\multicolumn{11}{l}{\scriptsize }
\\\hline
$C_4^{1*}$&1234&4231&3241&3214&1234&&&&&&
\\

$C_4^{2*}$&1234&4231&4321&1324&1234

\\

$C_4^{3*}$&1234&4231&3241&1243&1234
\\

$C_4^{4*}$&1234&4231&4213&1243&1234

\\ \hline

\end{tabular}

\vskip 0.2cm
Table 4.  $4$-cycles containing  $(1234,4231)$ in $BS_4$.
\end{center}

The cycles of every even length between $10$ and $24$ containing (1234,4231) can be constructed by the same argument as that of Subcase 1.1.\q

\vskip.2cm

{\bf Theorem 3.3} {\em For $n\ge 3$, the bubble-sort star graph $BS_n$ is
edge-bipancyclic.
Furthermore,
every edge of $BS_n$ lies on at least four different cycles of every even length $l$ with $4\leq l\leq n!$.}

\vskip.2cm

\p We prove this theorem by induction on $n$. For $n=3,4$, the result holds by Lemmas 3.1 and 3.2. Assume $n\geq5$.
Let $e=(u,v)\in E(BS_n)$ be an arbitrary edge. Since $BS_n$ is
vertex transitive, we can assume that $u=12\cdots n$. Then $u\in
V(BS_n(n))$. Now we will consider the following three cases.

{\bf Case 1.} {\em $v\in V(BS_n(n))$.}

Note that $BS_n(n)\cong BS_{n-1}$. Hence $e=(u,v)$
lies on at least four different cycles of every even length between $4$ and $(n-1)!$ in $BS_n(n)$ by  induction hypothesis. Let $q$ be an integer with $1\leq q\leq n-1$ and $p$ be an even integer with $2\leq p\leq (n-1)!$. We first prove the following claim.

{\bf Claim 1.} {\em There are at least four different cycles of length $l$ containing  $e$, where $l=q(n-1)!+p$.}

{\bf Proof of Claim 1.} We prove this claim by induction on $q$.

Suppose $q=1$.
Let $H_n$ be a Hamiltonian cycle in $BS_n(n)$ containing $e$. By Lemma 2.3, there exists an edge $e_1=(s_1,r_1)\in E(H_n)$ with $e_1\neq e$ that has a coupled pair-edge $e_1'=(s_1',r_1')$ such that $e_1'\in E(BS_n(m))$ for every $1\leq m\leq n-1$. Thus the edge set $E(H_n)\cup\{e_1',(s_1,s_1'),(r_1,r_1')\}-\{e\}$ forms a $((n-1)!+2)-$cycle containing $e$ in $BS_n$. Since $e=(u,v)$
lies on at least four different Hamiltonian cycles in $BS_n(n)$, there are at least four $((n-1)!+2)-$cycles containing $e$ in $BS_n$. Now we will consider the case that $4\leq p\leq (n-1)!$.

Let $e_{11}'=(s_{11}',r_{11}')\in E(BS_n(1))$ be a coupled pair-edge of  $e_{11}=(s_{11},r_{11})\in E(H_n)$ with $e_{11}\neq e$. Since $BS_n(1)\cong BS_{n-1}$, $e_{11}'$ can be contained in at least four different $p$-cycles in $BS_n(1)$, where $4\leq p\leq(n-1)!$. Let $C_{p}$ be any $p$-cycle in $BS_n(1)$ containing $e_{11}'$, then the edge set $E(H_n)\cup E(C_p)\cup\{(s_{11},s_{11}'),(r_{11},r_{11}')\}-\{e_{11}',e_{11}\}$ forms a $((n-1)!+p)$-cycle containing $e$. Thus $e$ lies on
at least four different cycles of length $(n-1)!+p$. Hence the claim holds for $q=1$.

Suppose the claim holds for $q=k-1$, where $2\leq k \leq n-1$. Thus there are at least four different cycles of every even length between $4$ and $k(n-1)!$ containing $e$.
Let $C_{k(n-1)!}$ be a cycle of length $k(n-1)!$ containing  $e$. Let $H_i$ be a Hamiltonian cycle in $BS_n(i)$ with $E(H_i)=(E(C_{k(n-1)!})\cap E(BS_n(i)))\cup\{e'_{i1},e_{(i+1)1}\}$, where  $e'_{i1}$ is a coupled pair-edge of $e_{i1}$ for every $i=1,2,\cdots,k-1$ (See Fig. 5).

\begin{figure}[ht]
   \begin{center}
    \includegraphics[scale=0.7]{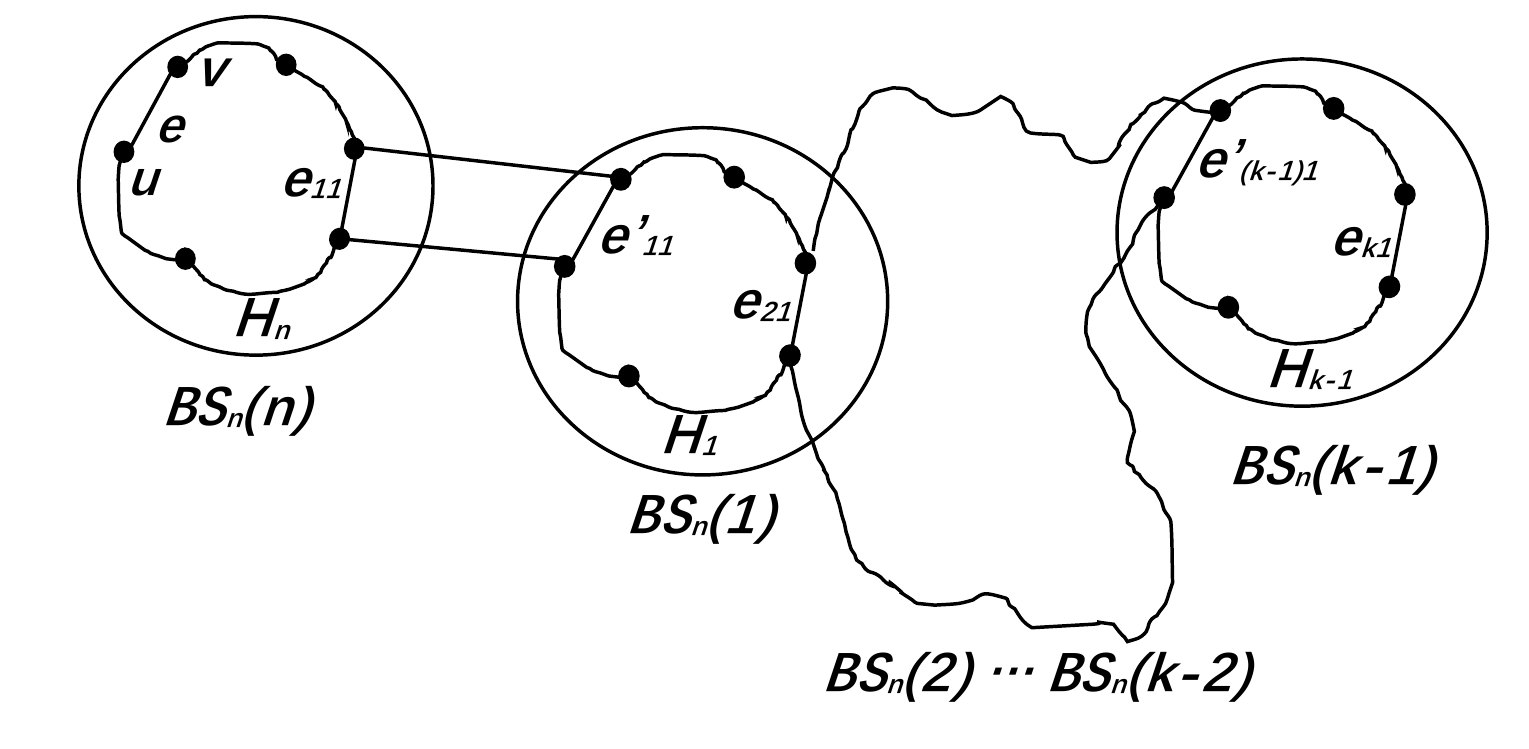}
    \end{center}
    \vspace{-0.5cm}\caption{\label{F2-5} Illustration of Claim 1.}
\end{figure}

Now we will prove the claim holds for $q=k$. By Lemma 2.3, there exists an edge $e*=(x,y)\in E(H_i)-\{e'_{i1},e_{(i+1)1}\}$ that has a coupled pair-edge $e'=(x',y')\in E(BS_n(j))$ for every $i\in\{1,2,\cdots,k-1,n\}$ and $j\in\{k,k+1,\cdots,n-1\}$. Thus the edge set $E(C_{k(n-1)!})\cup\{e',(x,x'),(y,y')\}-\{e*\}$ forms a $(k(n-1)!+2)$-cycle containing $e$. Since $k(n-k)\geq4$ for $n\geq5$ and $2\leq k \leq n-1$, there are at least four different $(k(n-1)!+2)-$cycles containing $e$ in $BS_n$. Now we will consider the case that $4\leq p\leq (n-1)!$.

Let $e_{k1}'=(s_{k1}',r_{k1}')\in E(BS_n(k))$ be a coupled pair-edge of  $e_{k1}=(s_{k1},r_{k1})\in E(H_{k-1})$, where $e_{k1}\neq e'_{(k-1)1}$. Since $BS_n(k)\cong BS_{n-1}$, by induction hypothesis, $e_{k1}'$ can be contained in at least four different $p$-cycles in $BS_n(k)$, where $4\leq p\leq(n-1)!$. Let $C_{p}$ be any $p$-cycle in $BS_n(k)$ containing $e_{k1}'$, then the edge set $E(C_{k(n-1)!})\cup E(C_p)\cup\{(s_{k1},s_{k1}'),(r_{k1},r_{k1}')\}-\{e_{k1}',e_{k1}\}$ forms a $(k(n-1)!+p)$-cycle containing $e$. Thus  $e$ lies on
at least four different cycles of length $k(n-1)!+p$. Hence the claim holds.\q

By Claim 1, we have that $e$ lies on at least four different cycles of every even length $l$ with $4\leq l\leq n!$.

{\bf Case 2.} {\em $v=u^-$.}

The four different $4$-cycles containing  $e$
in $BS_n$ are listed in Table 5.

\begin{center}
\setlength{\tabcolsep}{3.5pt}\renewcommand\arraystretch{0.8}
\begin{tabular}{cllllllllllll}\hline

\multicolumn{1}{c}{$cycle$}&\multicolumn{11}{l}{\scriptsize }
\\\hline
$C_4^{1}$&$~u$& $~u^-$&$~u^-\circ(1,2)$&$~u\circ(1,2)$&$~u$&&&&&&
\\

$C_4^{2}$&$~u$& $~u^-$&$~u^-\circ(1,3)$&$~u\circ(1,3)$&$~u$

\\

$C_4^{3}$&$~u$& $~u^-$&$~u^-\circ(2,3)$&$~u\circ(2,3)$&$~u$
\\

$C_4^{4}$&$~u$& $~u^-$&$~u^-\circ(1,n-1)$&$~u\circ(1,n-1)$&$~u$

\\ \hline

\end{tabular}

\vskip 0.2cm

Table 5.  $4$-cycles containing  $(u,u^-)$ in $BS_n$.
\end{center}

Note that $(u, u\circ(1,2))\in E(BS_n(n))$ is an edge of $C_4^{1}$ in Table 5. Since $BS_n(n)\cong BS_{n-1}$, by induction hypothesis, $(u, u\circ(1,2))$ can be contained in four different $p$-cycles in $BS_n(n)$ for $p=4,6,\cdots,(n-1)!$. Let $C^{p'}$ be any $p$-cycle in $BS_n(n)$ containing $(u, u\circ(1,2))$, then the edge set $E(C_4^{1})\cup E(C^{p'})-\{(u, u\circ(1,2))\}$ forms a $(2+p)$-cycle containing $e$. Thus $e$ lies on
at least four different cycles of every even length between 6 and $(n-1)!+2$.

The cycles of every even length between $(n-1)!+4$ and $n!$ containing $e$ can be constructed by the same argument as that of Case 1.

{\bf Case 3.} {\em $v=u^+$.}

The four different $4$-cycles containing $e$
in $BS_n$ are listed in Table 6.

\begin{center}
\setlength{\tabcolsep}{3.5pt}\renewcommand\arraystretch{0.8}
\begin{tabular}{cllllllllllll}\hline

\multicolumn{1}{c}{$cycle$}&\multicolumn{11}{l}{\scriptsize }
\\\hline
$C_4^{1*}$&$~u$& $~u^+$&$~u^+\circ(2,3)$&$~u\circ(2,3)$&$~u$&&&&&&
\\

$C_4^{2*}$&$~u$& $~u^+$&$~u^+\circ(3,4)$&$~u\circ(3,4)$&$~u$&&&&&&

\\

$C_4^{3*}$&$~u$& $~u^+$&$~u^+\circ(n-1,n)$&$~u\circ(n-1,n)$&$~u$&&&&&&
\\

$C_4^{4*}$&$~u$& $~u^+$&$~u^+\circ(n-1,n)$&$~u\circ(1,n-1)$&$~u$&&&&&&

\\ \hline

\end{tabular}

\vskip 0.2cm

Table 6.  $4$-cycles containing  $(u,u^+)$ in $BS_n$.
\end{center}

Note that $(u, u\circ(2,3))\in E(BS_n(n))$ is an edge of $C_4^{1*}$ in Table 6. Since $BS_n(n)\cong BS_{n-1}$, by induction hypothesis, $(u, u\circ(2,3))$ can be contained in four different $p$-cycles in $BS_n(n)$ for $p=4,6,\cdots,(n-1)!$. Let $C^{p''}$ be any $p$-cycle in $BS_n(n)$ containing $(u, u\circ(2,3))$, then the edge set $E(C_4^{1*})\cup E(C^{p''})-\{(u, u\circ(2,3))\}$ forms a $(2+p)$-cycle containing $e$. Thus  $e$ lies on
at least four different cycles of every even length between 6 and $(n-1)!+2$.

The cycles of every even length between $(n-1)!+4$ and $n!$ containing $e$ can be constructed by the same argument as that of Case 1.

Thus we complete the proof of our main theorem.\q
\vskip.2cm

By Theorem 3.3, we immediately obtain the following corollaries.

\vskip.2cm

{\bf Corollary 3.4} {\em For $n\ge 3$, the bubble-sort star graph $BS_n$ is
vertex-bipancyclic.}

\vskip.2cm
{\bf Corollary 3.5} {\em For $n\ge 3$, the bubble-sort star graph $BS_n$ is
bipancyclic.}
\vskip.3cm

\n{\large\bf Acknowledgements}
\vskip.2cm
This research is supported by National Natural Science Foundation of China (No. 11801450, 11771247), Natural Science Foundation of  Shaanxi Province, China (No. 2019JQ-506).

\end{document}